 \documentstyle[12pt,amscd]{amsart}
\newtheorem{prop}{Proposition}
\newtheorem{thm}{Theorem}
\newtheorem{lem}{Lemma}

\theoremstyle{definition}

\newtheorem{df}{Definition}

\theoremstyle{remark}

\newtheorem{rem}{Remark} 
\newcommand{\N}{{\Bbb N}}

\newcommand{\dint}{\displaystyle\int}

\makeatletter
\def\@currentlabel{2.1}\label{e:dispaa}
\def\@currentlabel{2.21}\label{e:dispau}
\def\@currentlabel{2.22}\label{e:dispav}
\def\@currentlabel{2.23}\label{e:dispaw}
\def\@currentlabel{2.24}\label{e:dispax}
\def\theequation{\@arabic\c@equation}

\makeatother

\makeatletter
\def\alphenumi{%
  \def\theenumi{\alph{enumi}}%
  \def\p@enumi{\theenumi}%
  \def\labelenumi{(\@alph\c@enumi)}}
\makeatother

\begin{document}

\title{COMPLEMENTED COPIES OF $\ell^1$ AND PELCZYNSKI'S PROPERTY
(V$^*$) IN  BOCHNER FUNCTION SPACES}
\author{Narcisse Randrianantoanina}
\address{Department of Mathematics, The University of Texas at  
Austin,
Austin, TX 78712}
\email{nrandri@@math.utexas.edu}
\maketitle

\vskip 3in
\begin{abstract}  Let $X$ be a Banach space and $(f_n)_n$ be a
bounded  sequence in $L^1(X)$. We prove a complemented
version of the celebrated Talagrand's dichotomy i.e
we show that if $(e_n)_n$
denotes the unit vector basis of $c_0$, there exists a sequence $g_n  
\in
\text{conv}(f_n,f_{n+1},\dots)$ such that for almost every $\omega$,  
either
the sequence $(g_n(\omega) \otimes e_n)$ is weakly Cauchy in $X
\widehat{\otimes}_\pi c_0$ or it is equivalent to the unit vector  
basis
of
$\ell^1$. We then get a criterion for a bounded sequence to contain  
a
subsequence equivalent to a  complemented copy of
$\ell^1$ in $L^1(X)$. As an application, we show that for a Banach  
space
$X$, the space $L^1(X)$ has
Pe\l czy\'nski's property $(V^*)$ if and only if $X$ does.
\end{abstract}

\section{INTRODUCTION}
Let $X$ be a Banach space and $(\Omega,\Sigma,\lambda)$ be a finite
measure space. If $1\leq p <\infty$, we denote by $L^p(\lambda,X)$
the Banach space of all (class of) $X$-valued $p$-Bochner integrable
functions with its usual norm. If $E$ and $F$ are Banach spaces, we  
denote by
$E \widehat{\otimes}_\pi F$ the projective tensor product of $E$ and  
$F$.
We will say that a sequence $(x_n)_n$ is equivalent to a  
complemented
 copy of $\ell^1$ if $(x_n)_n$ is equivalent to the unit vector  
basis
of $\ell^1$ and its closed linear span is complemented in $X$.

One of the many important problems in the theory of Banach spaces is
to recognize different structure of subspaces of a given space.
In this paper, we will be mainly conserned with sequences in the  
Bochner
space $L^1(\lambda,X)$ that are equivalent to a complemented copy of
$\ell^1$. Let us recall that
in \cite{T2}, Talagrand proved a fundamental theorem
characterizing weakly Cauchy sequences and sequences that are  
equivalent
to the unit vector basis of $\ell^1$  in the Bochner space
$L^1(\lambda,X)$, relating a given sequence $(f_n)_n$ to its values
$(f_n(\omega))_n$ in $X$. Our main goal is to provide a complemented
version of Talagrand's result. One way one might tackle this problem  
is
to consider for a given $(f_n)_n$ in $L^1(\lambda,X)$, the  
corresponding
 sequence $(f_n \otimes e_n)_n$ in
$L^1(\lambda,X) \widehat{\otimes}_\pi c_0$ (which can be viewed as
the Bochner space $L^1(\lambda, X \widehat{\otimes}_\pi c_0)$),  
where
 $(e_n)_n$ is the unit vector basis of $c_0$.
 The basic motivation behind this approach is the
 well known fact that a bounded sequence $(x_n)_n$ in a Banach
space
$X$ contains a subsequence
equivalent to a complemented copy of $\ell^1$ if and only if the
sequence $(x_n \otimes e_n)_n$ is not a weakly null sequence in $X
\widehat{\otimes}_\pi c_0$. Therefore the behavior
of the sequence  $(f_n\otimes e_n)_n$ will determine whether or not  
the
sequence $(f_n)_n$ has a subsequence equivalent to a complemented  
copy
of $\ell^1$.
As in \cite{T2}, we try to relate the sequence $(f_n \otimes e_n)_n$
to its values $(f_n(\omega) \otimes e_n)_n$ in $X  
\widehat{\otimes}_\pi
c_0$ to see how a particular structure of the space $X$ can be  
carried
 on to the Bochner space $L^1(\lambda,X)$.
The main result of this paper is the following extension of  
Talagrand'
s theorem: Let $(f_n)_n$ be a
bounded sequence in $L^1(\lambda,X)$, then there exists a sequence
$g_n \in \text{conv}(f_n,f_{n+1},\dots)$ such that
$(g_n(\omega) \otimes e_n)_n$ is either weakly Cauchy or  equivalent  
to
the unit vector basis of $\ell^1$. We follow a line of reasoning
similar to that of Talagrand (\cite{T2}). Our  main focus
is to carry out all the steps in such a way that the convex
combination is taken only on the sequence $(f_n)_n$ not on the
sequence $(f_n \otimes e_n)_n$.

In Section 3., we apply our main theorem for the study of property
(V$^*$) introduced by Pe\l czy\'nski in \cite{PL1}. The most notable
examples of Banach spaces that have property (V$^*$) are  
$L^1$-spaces
and it is a natural question to ask for what Banach spaces $X$ the
space $L^1(\lambda,X)$ has property (V$^*$). The most one could hope
for is that $L^1(\lambda,X)$ has property (V$^*$) if and only if $X$
does.
  This question was  studied by several authors.
 Partial results can be found in \cite{B1},
 \cite{EM4}, \cite{SS4} and more recently in \cite{LEU}. We present  
a
complete positive answer to this question (see Theorem 2 below).

 Our notation and terminology are standard and can be found in
 \cite{D1} and \cite{DU}.

\section{COMPLEMENTED VERSION OF TALAGRAND'S THEOREM.}

By way of motivation, let us begin with the following well known
proposition that justifies our approach. The word operator will  
always
mean linear bounded operator and $\cal L(X,Y)$ will stand for the  
Banach
space of all operators from $X$ to $Y$.
\begin{prop}
Let $X$ be a Banach space and $(x_n)_n$ be a bounded sequence in $X$  
that
is equivalent to the $\ell^1$ basis.
Then the following statements are equivalent:
\begin{itemize}
\item[(i)] The sequence $(x_n)_n$ is equivalent to a complemented  
copy of
$\ell^1$;
\item[(ii)] There exists an operator $T \in \cal L(X,\ell^1)$ such  
that
$\langle Tx_n, e_n \rangle \geq 1$ for all $n \in \N$.
\end{itemize}
\end{prop}

Using the fact that the space $\cal L(X,\ell^1)$ is the dual of
$X \widehat{\otimes}_\pi c_0$, condition (ii) of
Proposition~1 can be restated as:
$$ (*)\ \  \text{There exists}\ T \in (X
\widehat{\otimes}_\pi c_0)^*\ \text{ such that}
\ \langle T, x_n \otimes e_n
\rangle \geq 1 \ \text{for all}\  n \in \N.$$
 Now (*)  implies that the
sequence $(x_n \otimes e_n)_n$ is not a weakly null sequence. So the
problem of whether or not $(x_n)_n$ is equivalent to a
complemented copy of $\ell^1$ in $X$ is reduced to the study of
weak convergence of $(x_n \otimes e_n)_n$ in
$X \widehat{\otimes}_\pi c_0$.

The following result is our main criterion for determining if a  
given
sequence in a Bochner space has a subsequence that is equivalent to  
a
complemented copy  of $\ell^1$.

\begin{thm} Let $X$ be a Banach space and $(\Omega,
\Sigma,\lambda)$ be a probability space. Let $(f_n)_n$ be a bounded
sequence in $L^1(\lambda,X)$. Then there exist a sequence $g_n \in
\text{conv}(f_n, f_{n+1},\dots)$ and  two measurable subsets $C$ and
$L$ of $\Omega$ with $\lambda(C \cup L)=1$ such that:
\begin{itemize}
\item[(a)] for $\omega \in C$, the sequence $\left(g_n(\omega)  
\otimes
e_n \right)_n$ is weakly Cauchy in the space $X  
\widehat{\otimes}_\pi
c_0$.
\item[(b)] for $\omega \in L$, the sequence
$\left(g_n(\omega) \otimes
e_n \right)_n$   is equivalent to the $\ell^1$ basis in $X
\widehat{\otimes}_\pi c_0$.
\end{itemize}
\end{thm}

The proof uses many (if not all) ideas from Talagrand's theorem so
we recommend that the reader should get familiar to its proof first
before reading our extension.
However because of the complexity of the proof of
Talagrand's theorem, we decided to present all critical details.

Using similar argument as in \cite{DSC}, we can assume without loss  
of
generality that the sequence $(f_n)_n$ is such that $\sup\limits_{n  
\in
\N} ||f_n||_\infty \leq 1$.

For convenience, we will use the following notation:
\begin{itemize}
\item[(i)] For two sequences $(g_n)_n$ and $(f_n)_n$, we write
$(g_n) \ll (f_n)$ if there exists $k \in \N$ so that $\forall n \geq  
k$,
$g_n \in \text{conv}(f_n, f_{n+1},\dots)$ and by passing to a
subsequence (if necessary), we will always assume that there exist  
two
sequences of integers $(p_n)$ and $(q_n)_n$ such that $p_1\leq
q_1<p_2\leq q_2\dots$ and $g_n=\sum_{i=p_n}^{q_n} a_i f_i$;
\item[(ii)] For a Banach space $Y$, we denote by $Y_1$ the closed
unit ball of $Y$.
\end{itemize}

\noindent
Case 1: The Banach space $X$ is separable;

 If the space $X$ is separable, so is the space $X  
\widehat{\otimes}_\pi
c_0$
and therefore $\cal L(X,\ell^1)_1$ the unit ball of $\cal  
L(X,\ell^1)=
(X\widehat{\otimes}_\pi c_0)^*$
endowed with the weak$^*$-topology is compact metrizable.

Let us now consider $(U_n)_n$ a countable basis for the weak$^*$-  
topology
on $\cal L(X,\ell^1)_1$.

Following Talagrand \cite{T2}, we denote by $\cal K$ the set of all
(weak$^*$) compact sets of $\cal L(X,\ell^1)_1$ and we say that a  
map
$\omega \to K(\omega)$ from $\Omega$ to $\cal K$ is measurable if for  
each $n \in \N$, the set
$\left\{\omega \in \Omega;\ {K(\omega)} \cap {U_n}\neq \emptyset
 \right\}$
is measurable.

As in \cite{T2}, we will make use of the following elementary lemma:
\begin{lem}
If for each $k \in \N$, we have $(f_n^{k+1}) \ll (f_n^{k})$,
then there exists a sequence $(k_n)$ such that if we set
$g_n= f_n^{k_n}$ we have $(g_n) \ll (f_n^{k})$ for all $k \in \N$.
\end{lem}
We are now ready to begin the proof of the Theorem.

Let $\omega \to K(\omega)$ from $\Omega$ to $\cal K$ be a measurable  
map
and $V$ be a weak$^*$-open subset of $\cal L(X,\ell^1)_1$ and let
$g_n: \Omega \to X_1$ be  a bounded sequence in  
$L^\infty(\lambda,X)$
such that $(g_n) \ll (f_n)$. Let $g_n =\sum\limits_{i=p_n}^{q_n}
\lambda_i f_i $ be the representation of $g_n$ as block convex
combination of the $f_n$'s.

We set:
\begin{equation}
\overline{g_n}(\omega)=\sup\limits_{k \geq q_n}\sup\{\langle
T(g_n(\omega)),e_k
\rangle,\ T
\in
V \cap K(\omega)\}
\end{equation}
\begin{equation}
\theta(g)(\omega)=\limsup_{n \to \infty} \overline{g_n}(\omega).
\end{equation}

Notice that the definition of $\overline{g_n}$ depends on the
representation of $g_n$ as block convex combination of the $f_n$'s.
It is clear that $||\overline{g_n}||_\infty \leq 1$ and we claim  
that
$\overline{g_n}$ is measurable. To see this notice that for each
$k \in \N$, the map $\omega \to \sup\{ \langle T, g_n(\omega)\otimes
e_k \rangle,\ T \in V \cap K(\omega)\}$ was already proved to be
measurable by Talagrand so the claim follows.

\begin{lem}
 There exists $(g_n) \ll (f_n)$ such that  if
$(h_n) \ll (g_n)$ we have $\lim\limits_{n \to \infty} ||\theta(g)-
\overline{h_n} ||_1 =0$.
\end{lem}
\begin{pf}
The proof is done more or less the same as in \cite{T2}; let
$g_n^1 = f_n$ and construct by induction sequences $g^p =(g_n^p)$  
such
that for $p\geq 1$, one has
$$\int \theta(g^p)(\omega) \ d\lambda(\omega) \leq 2^{-p} +\inf\{
\int \theta(\phi)(\omega) \ d\lambda(\omega); \ \phi \ll  
g^{p-1}\}.$$
By Lemma 1., there is a sequence $(g_n)$ such that $g \ll g^p$ for
each $p \in N$; in particular $g \ll u$.

Let $h \ll g$ . We claim that $\theta(h)=
\theta(g)$. To see the claim write $h_n =\sum\limits_{i=p_n}^{q_n}
\alpha_i f_i$ ; $h_n= \sum\limits_{j=a_n}^{b_n} \beta_j g_j$ and
$g_n = \sum\limits_{l=c_n}^{d_n} \gamma_l f_l$. We have
\begin{align*}
\overline{h_n}(\omega)
 &=\sup\limits_{k \geq q_n}\sup\{\sum\limits_{i=a_n}^{b_n} \beta_i
\langle T(g_i(\omega)), e_k
\rangle,  \ T \in V \cap K(\omega)\} \\
 &\leq \sup\limits_{k \geq q_n}\sup\{\sup\limits_{i \in [a_n,b_n]}
\langle T(g_i(\omega)),e_k
 \rangle,\ T \in V \cap K(\omega)\}\\
\intertext{for each $n \in \N$ and $\omega \in \Omega$,
 so there exist $T \in V\cap K(\omega)$, $i_n \in [a_n,b_n]$
 such that}
 \overline{h_n}(\omega)
 &\leq \sup\limits_{k\geq q_n}\langle T(g_{i_n}(\omega)),e_k \rangle
+ \frac{1}{2^n}\\
\intertext{ but $q_n \geq d_{b_n} \geq d_{j_n}$ so we get that }
 \overline{h_n}(\omega)
&\leq \sup\limits_{k\geq d_{j_n}}\langle T(g_{i_n}(\omega)),e_k  
\rangle
+ \frac{1}{2^n}\\
 &\leq \overline{g_{i_n}}(\omega) + \frac{1}{2^n}
\end{align*}
 and by taking the limsup, we get that $\theta(h) \leq \theta(g)$.

 In the other hand we have  for each $p \in \N$,
\begin{align*}
\inf\{\int \theta(\phi)(\omega)\ d\lambda(\omega),\ \phi \ll  
g^{p-1}\}
&\leq \int \theta(h)(\omega)\ d\lambda(\omega) \\
&\leq \int\theta(g)(\omega)\ d\lambda(\omega)\\
&\leq \inf\{\int \theta(\phi)(\omega)\ d\lambda(\omega),\ \phi \ll
g^{p-1} \} +2^{-p}
\end{align*}
hence $\dint \theta(h)(\omega)\ d\lambda(\omega)=\dint  
\theta(g)(\omega)\
d\lambda(\omega)$.

\noindent We claim that $\theta(g)= \lim\limits_{n \to \infty}
\overline{h_n}$ for the weak$^*$-topology in $L^\infty(\lambda)$:  
for
that let $\phi$ be a cluster point of $(\overline{h_n})_n$. Since
$\theta(h)= \limsup\limits_{n \to \infty} \overline{h_n} \leq  
\theta(g)$,
one has $\phi \leq \theta(g)$. Moreover if we choose
${h_n}^\prime = \sum_{i=a_n}^{b_n} \alpha_i h_i $ such that $ a_1  
\leq b_1
<a_2 \leq b_2 < \dots$ and $||\sum_{i=a_n}^{b_n}  
\alpha_i\overline{h_i} -\phi
||_1 \leq 2^{-n}$, we have $\lim\limits_{n\to \infty}
\sum_{i=a_n}^{b_n} \alpha_i\overline{h_i}(\omega) =\phi(\omega)$\ a.e  
but
for any $n \in \N$, the above  estimate shows that:

\noindent  $\overline{h_n^\prime}(\omega)\leq
\overline{h_{i_n}}(\omega)
+2^{-n}$
and hence $\theta(g)=\theta(h^\prime)=\limsup\limits_{n \to \infty}
\overline{h_n^\prime} \leq \phi \leq \theta(g)$ which shows that
$\theta(g) =\phi$\ a.e and the claim is proved.

To conclude the proof of the lemma, notice that $$\limsup\limits_{n  
\to
\infty} \overline{h_n} \leq \theta(g)$$ and $\lim\limits_{n \to  
\infty}
\dint \overline{h_n}(\omega)\ d\lambda(\omega) = \dint  
\theta(g)(\omega)\
d\lambda(\omega)$ so we get that  $\lim\limits_{n \to
\infty}||\overline{h_n}-
\theta(g)||_1 =0$
 and the lemma is proved.
\end{pf}

In a similar fashion we set for $g_n (\omega) =\sum_{i=p_n}^{q_n}
\lambda_i f_i$ a block convex combination of $f_n$'s
\begin{equation}
\widetilde{g_n}(\omega)=\inf\limits_{k \ge q_n}\inf\{\langle
T(g_n(\omega)), e_k
\rangle,\ T \in V \cap K(\omega)\}
\end{equation}

\begin{equation}
\varphi(g)(\omega)=\liminf\limits_{n \to \infty}  
\widetilde{g_n}(\omega)
\end{equation}
We have the corresponding lemma:
\begin{lem}
 There exists $(g_n)\ll (f_n)$ such that if $(h_n) \ll
(g_n)$ we have $\lim\limits_{n \to \infty}||\varphi(g)-\widetilde{
h_n}||_1=0$.
\end{lem}

We are now ready to present the main construction of the proof.
 Let us fix $a <b$ and let $\tau$ be the first
 uncountable ordinal. Set  $h_n^0 =f_n$ and
 $ \ K_0(\omega)=\cal L(X,\ell^1
 )_1$. For $\alpha <\tau$, we will construct (as in \cite{T2})
 sequences $h^\alpha =(h_n^\alpha)$, and measurable maps
 $K_\alpha: \Omega \to \cal K$ with the following properties:
 \begin{equation}
 \text{for}\ \beta <\alpha <\tau, \ h^\alpha \ll h^\beta.
 \end{equation}
 For $\alpha <\tau$ and $h \ll f$ (say $h_n = \sum_{i=a_n}^{b_n}
\lambda_i f_i$ a representation of $(h_n)_n$ as a block convex
combination of $(f_n)_n$ ) \ we define:
 \begin{equation}
 \begin{split}
\overline{h}_{n,k,\alpha}(\omega)&=\sup\limits_{m \geq  
b_n}\sup\{\langle
T(h_n(\omega)),e_m
\rangle,\ T \in  U_k \cap K_\alpha(\omega)\}\\
\widetilde{h}_{n,k,\alpha}(\omega)&=\inf\limits_{m \geq  
b_n}\inf\{\langle
T(h_n(\omega)),e_m
 \rangle, \ T \in U_k \cap K_\alpha(\omega)\}\\
 \theta_{k,\alpha}(h)(\omega)&=\limsup\limits_{n \to \infty}
  \overline{h}_{n,k,\alpha}(\omega)\\
  \varphi_{k,\alpha}(h)(\omega)&=\liminf\limits_{n \to \infty}
  \widetilde{h}_{n,k,\alpha}(\omega)
  \end{split}
  \end{equation}
  then for each $\alpha$ of the form $\beta +1$ and each $h \ll  
h^\alpha$,
  we have $\lim\limits_{n \to \infty}||\theta_{k,\beta}(h^\alpha)-
  \overline{h}_{n,k,\beta} ||_1 =0$;\ $\lim\limits_{n \to \infty}
  ||\varphi_{k,\beta}(h^\alpha)-\widetilde{h}_{n,k,\beta}||_1=0$.

If $\alpha$ is limit, we set
\begin{equation}
 K_\alpha(\omega)=\bigcap_{\beta <\alpha}
K_\beta(\omega);
\end{equation}

If $\alpha=\beta +1$, we have
\begin{equation}
K_\alpha(\omega)=\{T \in K_\beta(\omega), \ T \in U_k \Rightarrow
\theta_{k,\beta}(h^\alpha)>b,\ \varphi_{k,\beta}(h^\alpha)<a\}.
\end{equation}
The construction is done by induction. Suppose that the construction
 has been done for each ordinal $\beta <\alpha$. If $\alpha$ is  
limit,
 we set $K_\alpha(\omega)=\bigcap_{\beta<\alpha}K_\beta(\omega)$.
 Let $\beta_n$ be an increasing sequence of ordinals with
  $\alpha=\sup\beta_n$. By Lemma 1, there exists $(h^\alpha)$ with
  $(h^\alpha)\ll (h^{\beta_n})$ for each $n \in \N$. Therefore for
  $\beta <\beta_n$, $(h^\alpha)\ll (h^{\beta_n})\ll (h^\beta)$ and
  hence $(h^\alpha)\ll (h^\beta)$ so (5) is satisfied.
  The construction is done in the case of limit ordinal.

  Suppose now that $\alpha=\beta +1$. Using Lemma 2. and Lemma 3.,
one can construct  a sequence $(g^k)$ with $(g^1)=(h^\beta)$,
$(g^{k+1}) \ll (g^k)$ and such that for $(h) \ll (g^k)$,
$\lim\limits_{n\to \infty}||\theta_{k,\beta}(g^k)-
\overline{h}_{n,k,\beta}||_1=0$; $\lim\limits_{n \to \infty}
||\varphi_{k,\beta}(g^k)-\widetilde{h}_{n,k,\beta}||_1=0$. Apply
Lemma 1. to get a sequence $(h^\alpha)$ with $(h^\alpha) \ll (g^k)$
for each $k \geq 1$ and we claim that $(h^\alpha)$ satisfy (6). To  
see
the claim let us fix $(h) \ll (h^\alpha)$. By the definition of
$(h^\alpha)$, we have for each $k\geq 1$, $(h)\ll(g^k)$. It follows
that $\lim\limits_{n\to \infty}||\theta_{k,\beta}(g^k)-
\overline{h}_{n,k,\beta}||_1=0$ and
$\lim\limits_{n\to \infty}||\varphi_{k,\beta}(g^k)-
\widetilde{h}_{n,k,\beta}||_1=0$. Since $(h^\alpha)\ll (g^k)$, we get  
that
 $\lim\limits_{n\to \infty}||\theta_{k,\beta}(g^k)-
\overline{h}_{n,k,\beta}^\alpha||_1=0$ and
$\lim\limits_{n\to \infty}||\varphi_{k,\beta}(g^k)-
\widetilde{h}_{n,k,\beta}^\alpha||_1=0$
which shows that $\theta_{k,\beta}(g^k)=\theta_{k,\beta}(h^\alpha)$  
and
$\varphi_{k,\beta}(g^k)=\varphi_{k,\beta}(h^\alpha)$
and the claim is proved.

Define now $K_\alpha(\omega)$ by (8). The measurability of
$K_\alpha(.)$ can be proved using similar argument as in \cite{T2}.
 The construction is complete.

\noindent Claim: There exists $\alpha <\tau$ such that
$K_\alpha(\omega)= K_{\alpha +1}(\omega)$ for a.e. $\omega \in  
\Omega$.

In fact if we set for each $k \geq 1$,
$\Omega_k^\alpha=\{\omega;\ U_k \cap K_\alpha(\omega)=\emptyset\}$
then for each $k\in \N$, the sequence
$(\lambda(\Omega_k^\alpha))_{\alpha <\tau}$ is increasing, hence
eventually constant. Fix $\alpha$ such that for each $k \in \N$,
we have $\lambda(\Omega_k^{\alpha +1})=\lambda(\Omega_k^\alpha)$.
It is clear that for $\omega \notin
\bigcup\limits_{k\geq 1} (\Omega_k^{\alpha +1}\setminus
\Omega_k^\alpha)$, we have $K_\alpha(\omega)=K_{\alpha +1}(\omega)$.
The claim is proved.

We now set $(h)=(h^{\alpha +1})$, $C=\{\omega;\ K_\alpha(\omega)
=\emptyset\}$ and $M=\{\omega;\  K_\alpha(\omega)=K_{\alpha  
+1}(\omega)
\neq \emptyset\}$. Clearly $\lambda(C \cup M)=1$ and for the rest of
the proof we set
$h_n=
\sum_{i=p_n}^{q_n} \lambda_i f_i$  be a representation of $(h_n)_n$  
as
block convex combunation of $(f_n)_n$. We have the following
property of the measurable subset $C$:
\begin{lem}
If $\omega \in C$ and $ T \in \cal L(X,\ell^1)_1$ and $u \ll h$,  
then
either
\begin{itemize}
\item[(a)] $\limsup\limits_{n \to \infty} \langle T(u_n(\omega)),e_n
\rangle \leq b$ or
\item[(b)]  $\liminf\limits_{n \to \infty} \langle  
T(u_n(\omega)),e_n
\rangle \geq a$.
\end{itemize}
\end{lem}
\begin{pf}
Let $\omega \in C$ and $T \in \cal L(X,\ell^1)_1$.
Fix $u \ll h \ (\ll f)$ say $u_n = \sum_{i=a_n}^{b_n} \alpha_i f_i$  
for
all $n \in \N$. Consider $S: c_0 \longrightarrow c_0$ defined as  
follows:
\ $Se_{b_n}=e_n \ \forall n \in \N$ and $Se_j=0$ for $j \neq b_n$.  
The
operator $S$ is trivially linear and $||S||=1$.
Since
$S^* \circ T \in K_0(\omega)$ and $S^*\circ T \notin  
K_\alpha(\omega)$,
there is
a least ordinal $\beta$ for which $S^*\circ T\notin K_\beta(\omega)$.  
The
ordinal
$\beta$ cannot be a limit so $\beta=\gamma +1$ and $S^*\circ T \in
K_\gamma(\omega)$. By the definition of $K_\beta(.)$, there exists
$k \in \N$ with $S^* \circ T \in U_k$ but either
$\theta_{k,\gamma}(h^\beta)(\omega) \leq b$ or
$\varphi_{k,\gamma}(h^\beta)(\omega)\geq a$. Now since
$u\ll h^\beta$, we get that either
$$\limsup_{n\to \infty} \langle T(u_n(\omega)),e_n \rangle
=\limsup_{n \to \infty}\langle S^*\circ T(u_n(\omega)),e_{b_n}
\rangle
\leq \theta_{k,\gamma}(u)(\omega)\leq
\theta_{k,\gamma}(h^\beta)(\omega) \leq b$$
or
$$\liminf_{n \to \infty}\langle T(u_n(\omega)), e_n \rangle =
\liminf_{n \to \infty}\langle S^*\circ T(u_n(\omega)), e_{b_n}\rangle  
\geq
\varphi_{k,\gamma}(u)(\omega)\geq \varphi_{k,\gamma}(h^\beta)\geq  
a.$$
The lemma is proved.
\end{pf}
For the set $M$, we have the following lemma:
\begin{lem}
There exists a subsequence $(n(i))$ of the integers so that
  for almost every $\omega \in M$, there exists $k \in \N$ such
that the sequence $(h_{n(i)}(\omega) \otimes e_i)_{i \geq k}$ is  
$\delta$
equivalent to the $\ell^1$-basis in $X \widehat{\otimes}_\pi c_0$,
where $\delta=(b-a)/2$.
\end{lem}
\begin{pf}
Again we adopt the methods in \cite{T2} to our situation. Let us
denote by $F$ the set of finite sequences of zeroes and ones. For
$s \in F$, we will denote by $|s|$ the length of $s$. For
$s=(s_1,\dots,s_n)$ and $r=(r_1,\dots,r_m)$ with $n \leq m$, we say
that $s<r$ if $s_i=r_i$ for $i\leq n$. We will construct two  
sequences
of  integers $n(i)$, $m(i)$, measurable sets $B_i \subset M$ and
 measurable maps $Q(s,.): M \to \N$ such that the following
conditions are satisfied:
\begin{equation}
q_{n(1)}<m(1)<q_{n(2)}<m(2)<\dots <m(i)<q_{n(i+1)}<\dots
\end{equation}
\begin{equation}
\forall s \in F, \quad \sup\{Q(s,\omega);\ \omega \in M\} <\infty;
\end{equation}
\begin{equation}
\lambda(M\setminus B_i) \leq 2^{-i};
\end{equation}
\begin{equation}
\text{For}\ s,\ r\in F,\ s <r,\ \text{and}\
\omega \in \bigcap_{|s|\leq i \leq |r|} B_i, \ \text{one has}
\quad  U_{Q(r,\omega)}
\subset U_{Q(s,\omega)};
\end{equation}
\begin{equation}
\forall \omega \in M,\ s \in F,\quad
K_\alpha(\omega)\cap U_{Q(s,\omega)} \neq \emptyset;
\end{equation}
\begin{equation}
\begin{split}
&\forall p,\  \forall i  \leq p,\  \forall\  \omega \in
\bigcap_{i \leq j \leq p}B_j, \\
s_i=1 &\Rightarrow \forall T \in U_{Q(s,\omega)},\
 \sup\limits_{q_{n(i)} \leq k \leq m(i)} \langle
T(h_{n(i)}(\omega)),e_k\rangle
\geq b \\
 s_i=0 &\Rightarrow \forall T \in U_{Q(s,\omega)},\
 \inf\limits_{q_{n(i)}\leq k \leq m(i)} \langle  
T(h_{n(i)}(\omega)),e_k
\rangle
\leq a.
\end{split}
\end{equation}

Again the construction is done by induction. Before doing so we need  
the
following notation:

Let $n \in \N$, $j \in \N$ and $\alpha <\tau$.
Fix $m \geq n$, the following notation will be used.
\begin{align*}
  \overline{h}_{n,j,\alpha}^{(m)}(\omega)&=\sup\limits_{q_n\leq
k\leq m} \sup\{\langle T(h_n(\omega)),e_k \rangle,\ T \in
U_j \cap K_\alpha(\omega)\}\\
\widetilde{h}_{n,j,\alpha}^{(m)}(\omega)&=\inf\limits_{q_n\leq k  
\leq
m} \inf\{\langle T(h_n(\omega)),e_k \rangle,\ T \in
U_j \cap K_\alpha(\omega)\};
\end{align*}
It is clear that
$\overline{h}_{n,j,\alpha}(\omega)=\lim\limits_{m \to \infty}
\overline{h}_{n,j,\alpha}^{(m)}(\omega)$ a.e and
$\widetilde{h}_{n,j,\alpha}(\omega)=\lim\limits_{m \to \infty}
\widetilde{h}_{n,j,\alpha}^{(m)}(\omega)$ a.e .

\noindent
For $i=1$, recall that
$U_0=\cal L(X,\ell^1)_1$. Since $K_{\alpha +1}(\omega) \neq  
\emptyset$
for $\omega \in M$, one has $\theta_{0,\alpha}(h)(\omega)>b$ and
$\varphi_{0,\alpha}(h)(\omega)<a$ but since
$$\lim_{n \to
\infty}||\theta_{0,\alpha}(h)-\overline{h}_{n,0,\alpha}||_1=
\lim_{n \to
\infty}||\varphi_{0,\alpha}(h)-\widetilde{h}_{n,0,\alpha}||_1=0,$$
there exists an integer $n(1)$ such that if we set
$$B_1^"=\left\{\omega \in M,\ \overline{h}_{n(1),0,\alpha}(\omega)  
>b;\
\widetilde{h}_{n(1),0,\alpha}(\omega)<a \right\}$$
we have $\lambda(M \setminus B_1^")\leq 2^{-3}$.
Since
$\overline{h}_{n(1),0,\alpha}(\omega)=\lim\limits_{m \to \infty}
\overline{h}_{n(1),0,\alpha}^{(m)}(\omega)$ a.e and
$\widetilde{h}_{n(1),0,\alpha}(\omega)=\lim\limits_{m \to \infty}
\widetilde{h}_{n(1),0,\alpha}^{(m)}(\omega)$ a.e, there exists an
integer $m(1)>q_{n(1)}$ such that if we set
$$B_1^\prime=\left\{\omega \in M,\
\overline{h}_{n(1),0,\alpha}^{m(1)}(\omega)
>b;\
\widetilde{h}_{n(1),0,\alpha}^{m(1)}(\omega)<a \right\},$$
we have $\lambda(M \setminus B_1^\prime)\leq 2^{-2}$.

Now for $\omega \in B_1^\prime$, we have:
\begin{align*}
\sup\limits_{q_{n(1)}\leq k \leq m(1)}\sup\{\langle
T(h_{n(1)}(\omega)), e_k \rangle,\ T \in K_\alpha(\omega)\} &>b \\
\intertext{and}
\inf\limits_{q_{n(1)} \leq k \leq m(1)}\inf\{\langle
T(h_{n(1)}(\omega)), e_k \rangle,\ T \in K_\alpha(\omega)\} &<a.
\end{align*}

For each $x \in X$, the maps  $T \to \sup\limits_{q_{n(1)} \leq k  
\leq
m(1)}\langle
Tx,e_k \rangle$  and   $T \to \inf\limits_{q_{n(1)}\leq k  \leq
m(1)}\langle Tx,
e_k \rangle$ are  continuous so
the sets
\begin{align*}
 \{T \in \cal L(X, \ell^1)_1,\ &\sup\limits_{q_{n(1)} \leq k \leq  
m(1)}
\langle Tx, e_k \rangle >b \} \\
\{T \in \cal L(X, \ell^1)_1,\ &\inf\limits_{q_{n(1)}\leq k \leq  
m(1)}
\langle Tx, e_k \rangle <a \}
\end{align*}
are open subsets. By a standard techniques
one can choose measurable maps $Q_0(.)$ and $Q_1(.)$ from $M$ to
$\N$ such that
\begin{align*}
T \in U_{Q_1(\omega)} &\Rightarrow \sup_{q_{n(1)}\leq k \leq m(1)}
\langle T(h_{n(1)}(\omega)), e_k \rangle >b;\ U_{Q_1(\omega)} \cap
K_\alpha(\omega) \neq \emptyset \\
T \in U_{Q_0(\omega)} &\Rightarrow \inf_{q_{n(1)}\leq k \leq m(1)}
\langle T(h_{n(1)}(\omega)), e_k \rangle <a;\ U_{Q_0(\omega)} \cap
K_\alpha(\omega) \neq \emptyset.
\end{align*}
There  exists an integer $l$ such that if $B_1 =\{\omega \in  
B_1^\prime;\
Q_0(\omega)<l,\ Q_1(\omega)<l\}$ we have $\lambda(M\setminus B_1)  
\leq
2^{-1}$. We define $Q((0),\omega)=Q_0(\omega)$ and
$Q((1),\omega)=Q_1(\omega)$ for $\omega \in B_1$ and
$Q((0),\omega)=Q((1), \omega)=0$ for $\omega \in M\setminus B_1$.  
The
required conditions (9)-(14) are satisfied.

 Suppose now that the result has been proved for $i$.
Let $l =\sup\{Q(s,\omega),\ |s|=i,\ \omega \in B_i\}$. Since
$K_\alpha(\omega)=K_{\alpha +1}(\omega)$, for $\omega \in M$,
condition (8)
implies that for each $k \in \N$, $$U_k \cap K_\alpha(\omega) \neq
\emptyset \Rightarrow \theta_{k,\alpha}(h)(\omega)>b,\
\varphi_{k,\alpha}(h)(\omega)<a.$$
We deduce as in the case $i=1$ that there is an integer $n(i+1)$  
such
that $q_{n(i+1)} > m(i)$ and the set

$$B_{i+1}^" =\left\{\omega \in M, \ \forall k \leq l, \ U_k \cap
K_\alpha(\omega)\neq \emptyset \Rightarrow
\overline{h}_{n(i+1),k,\alpha}(\omega)>b,\
\widetilde{h}_{n(i+1),k,\alpha}(\omega) <a \right\}$$

\noindent satisfies $\lambda(M \setminus B_{i+1}^") \leq 2^{-i -3}$.
Using similar argument as in the case $i=1$, one can pick an integer
$m(i+1)>q_{n(i+1)}$ so that the set

$$B_{i+1}^\prime =\left\{\omega \in B_{i+1}^", \ \forall k \leq l, \  
U_k
\cap K_\alpha(\omega)\neq \emptyset \Rightarrow
\overline{h}_{n(i+1),k,\alpha}^{(m(i+1))}(\omega)>b,\
\widetilde{h}_{n(i+1),k,\alpha}^{(m(i+1))}(\omega) <a \right\}$$
satisfies
$\lambda(M \setminus B_{i+1}^\prime)\leq 2^{-i-2}$.

For $\omega \in B_{i+1}^\prime$, one has in particular for $s \in  
F$,
$|s|=i$:
$$\overline{h}_{n(i+1),Q(s,\omega),\alpha}^{(m(i+1))}(\omega)  
>b;\quad
\widetilde{h}_{n(i+1),Q(s,\omega),\alpha}^{(m(i+1))}(\omega) <a.$$
It follows that for $s \in F$, $|s|=i$, there exist measurable maps
$Q_0(s,.)$ and $Q_1(s,.)$ from $M$ to $ \N$ suth that
\begin{align*}
T \in U_{Q_0(s,\omega)} &\Rightarrow \inf\limits_{q_{n(i+1)}\leq k  
\leq
m(i+1)}
\langle T(h_{n(i+1)}(\omega)), e_k \rangle <a, \\
 &U_{Q_0(s,\omega)} \cap K_\alpha(\omega) \neq \emptyset;\
U_{Q_0(s,\omega)} \subset U_{Q(s,\omega)} \\
\intertext{and}
T \in U_{Q_1(s,\omega)} &\Rightarrow \sup\limits_{q_{n(i+1)}\leq k  
\leq
m(i+1)}
\langle T(h_{n(i+1)}(\omega)), e_k \rangle >b, \\
 &U_{Q_1(s,\omega)} \cap K_\alpha(\omega) \neq \emptyset;\
U_{Q_1(s,\omega)} \subset U_{Q(s,\omega)}.
\end{align*}
There exists  an integer $l^\prime$ such that if we let
$$B_{i+1}= \left\{\omega \in B_{i+1}^\prime;\ \forall \ s \in F,\  
|s|=i,\
Q_0(s,\omega),\ Q_1(s,\omega)\leq l^\prime \right\}$$
then $\lambda(M\setminus
B_{i+1})\leq 2^{-i-1}$. The construction is done by setting
\begin{align*}
Q((s,0),\omega)&=Q((s,1),\omega)=0 \quad \text{if}\ \omega \in  
M\setminus
B_{i+1} \\
Q((s,0),\omega)&=Q_0(s,\omega);\ \  
Q((s,1),\omega)=Q_1(s,\omega)\quad
\text{if}\ \omega \in B_{i+1}.
\end{align*}

Let $L= \bigcup\limits_{k}\bigcap\limits_{i\geq k} B_{i}$. It is  
clear
that
$\lambda(M\setminus L)=0$ and we claim that if $\omega \in
\bigcap\limits_{i\geq k}B_i$, the sequence $(h_{n(i)}(\omega)\otimes
e_i)_{i\geq k}$ is $\delta$-equivalent to the unit vector basis of
$\ell^1$ in the Banach space $X \widehat{\otimes}_\pi c_0$. To see  
the
claim, let
$p\geq k$ and consider
 a subset $P$ of
$[k,p]$. Let $s \in F$ be a  sequence with $|s|=p$ and satisfies  
$s_i=1$
if $i \in P$, $s_i=0$ if $i\notin P$. From (14), there exists $T \in  
\cal
L(X,\ell^1)_1$ with:
\begin{align*}
k \leq i \leq p,\quad i \in P &\Rightarrow \sup\limits_{q_{n(i)}\leq  
m
\leq m(i)}\langle T(h_{n(i)}(\omega)), e_m \rangle \geq b \\
k \leq i \leq p,\quad i \notin P &\Rightarrow  
\inf\limits_{q_{n(i)}\leq m
\leq m(i)}\langle T(h_{n(i)}(\omega)), e_m \rangle  \leq a.
\end{align*}
Now for $i \in [k,p]$, choose $k(i)\in [q_{n(i)},m(i)]$ so that
\begin{align*}
k \leq i \leq p,\quad i \in P &\Rightarrow \sup\limits_{q_{n(i)}\leq  
m
\leq m(i)}\langle T(h_{n(i)}(\omega)), e_m \rangle= \langle
T(h_{n(i)}(\omega)), e_{k(i)}\rangle \\
k \leq i \leq p,\quad i \notin P &\Rightarrow  
\inf\limits_{q_{n(i)}\leq m
\leq m(i)}\langle T(h_{n(i)}(\omega)), e_m \rangle  =\langle
T(h_{n(i)}(\omega)), e_{k(i)}\rangle.
\end{align*}
By (9), the sequence $k(i)$ is increasing so there exists an  
operator
$S: c_0 \longrightarrow c_0$ of norm one such that $Se_i = e_{k(i)}$
 and it is now clear that:
\begin{equation}
\begin{split}
&k \leq i \leq p,\quad i \in P \Rightarrow
\langle S^* \circ T(h_{n(i)}(\omega)), e_i \rangle \geq b \\
&k \leq i \leq p,\quad i \notin P \Rightarrow
\langle S^* \circ T(h_{n(i)}(\omega)), e_i \rangle  \leq a.
\end{split}
\end{equation}
 And the claim follows from Rosenthal's argument in \cite{R1} (see  
also
\cite{D1} P.205). The lemma is proved.
 \end{pf}
 \begin{rem} Let $u\ll (h_{n(i)})_{i \in \N}$. Using the same  
argument
 as above, one can show that there exists a subsequence $(v_i)_i$ of
 $(u_i)_i$ such that $(v_i(\omega)\otimes e_i)_i$ is equivalent to  
the
 $\ell^1$ basis in $X \widehat{\otimes}c_0$ for a.e $\omega \in L$.
 \end{rem}
To complete the proof of the theorem, let $(a(k),b(k))$ be an  
enumeration
of all pairs of rational numbers with $a <b$. By induction we  
construct
sequences $g^k$ and measurable sets $C_k$, $L_k$ satisfying the
following:
\begin{itemize}
\item[(i)] $C_{k+1} \subset C_k, \ L_k \subset L_{k+1}$ and  
$\lambda(C_k
\cup L_k)=1$;
\item[(ii)]  $\forall \ \omega \in C_k$,\ $\forall m \leq k$, and $T  
\in
\cal L(X,\ell^1)_1$ then either
$\limsup\limits_{n \to \infty} \langle T(g_n^m(\omega)),e_n \rangle  
\leq
b(k)$ or
$\liminf\limits_{n \to \infty} \langle T(g_n^m(\omega)),e_n \rangle  
\geq
a(k);$
\item[(iii)] $\forall \ \omega \in L_m\setminus L_{m-1}$ with $2\leq  
m
\leq k$, the sequence $(g_n^k(\omega) \otimes e_n)_n $ is
$(b(m)-a(m))/2$-equivalent to the unit vector basis of $\ell^1$;
\item[(iv)] $g^{k+1} \ll g^k$.
\end{itemize}
Let $g^0 =f$, the steps above shows that one can find $g^1 \ll f$,
 measurable subsets $C_1$ and $L_1$
satisfying (i) - (iv). Suppose that $g^k$, $C_k$, $L_k$
have been constructed.
Again by the same reasoning for $a=a(k+1)$, $b=b(k+1)$ and
 $f_n=g_n^k$, there exist $g^{k+1} \ll g^k$ and measurable subsets
 $C_{k+1}$, $L_{k+1}$ with $\lambda(C_{k+1}\cup L_{k+1})=1$ and by
 Lemma 4. and Lemma 5., conditions (i) - (iv) are satisfied.

 We set $C=\bigcap\limits_{k\geq1} C_k$,
  $L=\bigcup\limits_{k\geq 1}L_k$ and $g_n= g_n^n$. It is clear that
$\lambda(C\cup L)=1$ and $g_n \ll g_n^k$ for each $k \in \N$;
in particular $g_n \ll f_n$. For $\omega \in C$, we have
 $(g_n(\omega) \otimes e_n)_n$ is weakly Cauchy. In fact since
 $g_n \ll g_n^k$, Lemma 4 asserts that for each
 $T \in \cal L(X,\ell^1)$ either
 $\limsup\limits_{n \to \infty}\langle T(g_n(\omega)),e_n \rangle
 \leq b(k)$ or
 $\liminf\limits_{n \to \infty}\langle T(g_n(\omega)),e_n \rangle
 \geq a(k)$ for all $k \in \N$. Hence
 $\limsup\limits_{n \to \infty} \langle T(g_n(\omega)),e_n \rangle=
 \liminf\limits_{n \to \infty}\langle T(g_n(\omega)),e_n \rangle$.
 Now for $\omega \in L$, there exists $k$ such that $\omega \in L_k$
 and since $g_n \ll g_n^k$, by Lemma 5,
 $(g_n(\omega)\otimes e_n)_{n\geq m}$ is equivalent to the unit
 vector basis of $\ell^1$ for some $m \in \N$. The proof
 of Theorem 1 is complete for the separable case.

\noindent
Case 2: General case;

One can reduce the general case to the separable one using the  
following
result of Heinrich
and Mankiewicz (see Proposition 3.4 of \cite{HM}):
\begin{lem} Let $X$ be a Banach space and $X_0$ be a separable  
subspace of
$X$. Then there exist a separable subpace $Z$ of $X$ that contains
$X_0$ and an isometric embeding $J: Z^* \to X^*$ such that
$\langle z,Jz^*\rangle=\langle z,z^*\rangle$ for every $z \in Z$ and
$z^* \in Z^*$. In particular $J(Z^*)$ is 1-complemented
in $X^*$.
\end{lem}

Let $(f_n)_{n \in \N}$ be a bounded sequence in $L^1(\lambda,X)$.  
Since
each $f_n$ has (essentially) separable range, there exists a  
separable
subspace $X_0$ of $X$ such that for a.e $\omega \in \Omega$,  
$f_n(\omega)
\in X_0$. Let $Z$ be a separable subspace as in the above lemma.   
The
sequence $(f_n)_n$ is bounded in $L^1(\lambda,Z)$ so by case 1,  
there
exist $g_n \in \text{conv}(f_n, f_{n+1},\dots)$, measurable subsets  
$C$
and $L$ of $\Omega$ with $\lambda(C \cup L)=1$ such that  for $\omega  
\in
C$, the sequence $(g_n(\omega) \otimes e_n)_n$ is weakly Cauchy in  
$Z
\widehat{\otimes}_\pi c_0$ and for $\omega \in L$, the sequence
$(g_n(\omega)\otimes e_n)_n$ is equivalent to the $\ell^1$ basis in  
$Z
\widehat{\otimes}_\pi c_0$.

We claim that the same conclusion holds if we replace $Z
\widehat{\otimes}_\pi c_0$ by $X \widehat{\otimes}_\pi c_0$. In fact  
if
$\omega \in C$ and $T \in \cal L(X,\ell^1)=(X \widehat{\otimes}_\pi
c_0)^*$,  the operator $T|_Z$ (the restriction of $T$ on $Z$) belongs  
to
$(Z \widehat{\otimes}_\pi c_0)^*$ so we have
\begin{align*}
\lim\limits_{n \to \infty}\langle T, g_n(\omega) \otimes e_n \rangle
&= \lim\limits_{n \to \infty} \langle T(g_n(\omega)),e_n \rangle \\
&=\lim\limits_{n \to \infty} \langle T|_Z, g_n(\omega) \otimes e_n
\rangle
\end{align*}
Hence $\lim\limits_{n \to \infty} \langle T, g_n(\omega) \otimes e_n
\rangle$ exists so the sequence $(g_n(\omega) \otimes e_n)_n$ is  
weakly
Cauchy in $X \widehat{\otimes}_\pi c_0$.

\noindent Now for $\omega \in L$, let $(a_n)_n$ be a finite sequence  
of
scalars. We have:
\begin{align*}
||\sum a_n g_n(\omega) \otimes e_n||_{X \widehat{\otimes}_\pi c_0}
&=\sup\{\sum a_n \langle T, g_n(\omega) \otimes e_n \rangle;\ T \in  
\cal
L(X,\ell^1)_1\} \\
&=\sup\{\sum a_n \langle  g_n(\omega), S e_n \rangle;\ S \in \cal
L(c_0,X^*)_1\} \\
&\geq \sup\{\sum a_n \langle  g_n(\omega), J\circ L e_n \rangle;\ L  
\in
\cal L(c_0,Z^*)_1\} \\
&= \sup\{\sum a_n \langle  g_n(\omega), L e_n \rangle;\ L \in
\cal L(c_0,Z^*)_1\} \\
&= ||\sum a_n g_n(\omega) \otimes e_n ||_{Z \widehat{\otimes}_\pi  
c_0}
\geq \delta  \sum |a_n|
\end{align*}
for some $\delta >0$. So the sequence $(g_n(\omega) \otimes e_n)_n$  
is
equivalent to the $\ell^1$ basis  in $X \widehat{\otimes}_\pi c_0$.
The theorem is proved.
\qed

 \noindent{\bf Remark:}
 In \cite{T2}, Talagrand extended his main theorem to the case
 of functions that are weak$^*$-scalarly measurable. It is not clear
 to us if one can get a similar result as in Theorem~1 for
 weak$^*$-scalarly measurable functions.

\section{APPLICATIONS: PROPERTY (V$^*$) AND (V$^*$)-SETS  FOR
$L^1(\lambda,X)$}

\begin{df} Let $X$ be a Banach space. A series $\sum_{n=1}^\infty  
x_n$
in $X$ is said to be a Weakly Unconditionally Cauchy (W.U.C.) if for  
every
$x^* \in X^*$, the series $\sum_{n=1}^\infty |x^*(x_n)|$ is  
convergent.
\end{df}

There are many criteria for a series to be a W.U.C. series (see for  
instance
\cite{D1} or \cite{WO}).

\begin{df} Assume that $X$ and $Y$ are Banach spaces. A bounded  
linear map
$T: X \to Y$ is said to be Unconditionally converging if $T$ sends
W.U.C. series in $X$ to unconditionally convergent series in $Y$.
\end{df}

In his fundamental paper \cite{PL1}, Pe\l czy\'nski proved the  
following
proposition:
\begin{prop} For a Banach space $X$, the following assertions are
equivalent:
\begin{itemize}
\item[(i)] A subset $H \subset X^*$ is relatively weakly compact  
whenever
$\lim\limits_{n \to \infty}\sup\limits_{x^* \in H} |x^*(x_n)|=0$ for
every W.U.C. series $\sum_{n=1}^\infty x_n$ in $X$ ;
\item[(ii)] For any Banach space $Y$, every bounded operator
$T: X \to Y$ that is unconditionally converging is weakly compact.
\end{itemize}
\end{prop}

\begin{df} A Banach space $X$ is said to have property (V) if it
satisfies one of the equivalent conditions of Proposition 3.
\end{df}

As a dual property, we have the following definition:
\begin{df} A Banach space $X$ is said to have property (V$^*$) if a  
subset $K$
of $X$ is relatively weakly compact whenever
$\lim\limits_{n \to \infty}\sup\limits_{x \in K} |x(x_n^*)|=0$ for
every W.U.C. series $\sum_{n=1}^\infty x_n^*$ in $X^*$.
\end{df}

\begin{df} A subset $K$ of a Banach space $X$ is called a  
(V$^*$)-set
if for every W.U.C. series $\sum_{n=1}^\infty x_n^*$ in $X^*$,
the following holds: $\lim\limits_{n\to \infty}\sup\limits_{x \in
K}|x(x_n^*)|=0$.
\end{df}

Hence a Banach space $X$ has property (V$^*$) if and only if every
(V$^*$)-set
 in $X$ is relatively weakly compact.

>From a result of Emmanuele \cite{EM4} (see also Godefroy and Saab
\cite{GOS1}), one can deduce
the following characterization of spaces that have property (V$^*$).

\begin{prop} A Banach space $X$ has property (V$^*$) if and only if
$X$ is weakly sequentially complete and given any sequence $(x_n)_n$
in $X$ that is equivalent to the unit vector basis of $\ell^1$,
there exists an operator $T: X \to \ell^1$ such that $(Tx_n)_n$ is  
not
relatively compact in $\ell^1$.
\end{prop}

The above proposition shows in particular that a Banach space $X$
has property (V$^*$) if and only if $X$ is weakly sequentially  
complete
and every sequence that is equivalent to the unit vector basis of
$\ell^1$ has a subsequence equivalent to a complemented copy of  
$\ell^1$.

In this section we will concentrate on property (V$^*$) and we shall  
refer
the reader to \cite{CKSS} and \cite{PL1} for more on property (V).

In \cite{SS4}, Saab and Saab showed (see Proposition 3. of  
\cite{SS4})
that a Banach space with the separable complementation property has
property (V$^*$) if and only if each  of its separable subspaces
has property
(V$^*$). On the next proposition, we will show that property (V$^*$)  
is
in fact separably determined.

\begin{prop} A Banach space $X$ has
property (V$^*$) if and only if all of its separable subspace has
property (V$^*$).
\end{prop}

\begin{pf} Since property (V$^*$) is easily seen to be stable by
subspaces,
 one implication is immediate.

 For the converse, we will use the  result of Heinrich
and Mankiewicz stated in Lemma 6 above.

Assume that every separable subspace of $X$ has property (V$^*$).
The space $X$ is trivially weakly sequentially complete. Let $K$ be
a bounded subset of $X$ that is not relatively weakly compact. There
exists a sequence $(x_n)_n$ in $K$ that is equivalent to the unit
vector basis of $\ell^1$ in $X$ and let
$X_0=\overline{\text{span}}\{x_n;\ n\in
\N\}$. The space $X_0$ is separable and consider $Z$ as in Lemma 6.
Since
$Z$ is separable, by assumption it has property (V$^*$) and
therefore there exists a W.U.C. series $\sum_{k=1}^\infty z_k^*$ in
$Z^*$ such that
$\limsup\limits_{k \to \infty}\sup\limits_{n \in \N} \langle  
z_k^*,x_n
\rangle >0$. Let $x_k^* =J(z_k^*)$; the series
$\sum_{k=1}^\infty x_k^*$ is a W.U.C. series in $X^*$ and
$$\langle x_k^*,x_n \rangle =\langle J(z_k^*),x_n \rangle=
\langle z_k^*,x_n \rangle.$$
So $\limsup\limits_{k\to \infty}\sup\limits_{n \in\N}\langle  
x_k^*,x_n \rangle>0$
which shows that $K$ is not a (V$^*$)-set.
\end{pf}

We are now ready to present the  main theorem of this section.

 Let  $E$ be a Banach
lattice with weak unit. By the classical representation (see  
\cite{LT}),
there exists a probability space $(\Omega,\Sigma,\lambda)$ such that
$L^\infty(\lambda) \subset E \subset L^1(\lambda)$ with $E$ being an
ideal and the inclusion being continuous.

Define  $E(X)$ to be the space of (class of) measurable map $f:  
\Omega \to X$
so that the measurable function $V(f)$ defined by
$V(f)(\omega)= ||f(\omega)||_{X}$ belongs to $E$. The space $E(X)$
endowed with the norm $||f||=||V(f)||_E$ is a Banach space.

 We have
the following stability result:

\begin{thm}
Let $X$ be a Banach space and $E$ be a Banach lattice that does not
contain any copy of $c_0$.
 The space $X$ has property (V$^*$) if and only if
$E(X)$ has property (V$^*$).
\end{thm}

\begin{pf} If $E(X)$ has property (V$^*$), then the space
$X$ has property (V$^*$) since property (V$^*$) is stable by  
subspace.

Conversely, assume that $X$ has property (V$^*$).
By Proposition 5., we
can assume without loss of generalities that
$E$ and $X$ are separable. By the classical representation, there
exists a probability space $(\Omega,\Sigma, \lambda)$ such that
$L^\infty(\lambda) \subset E \subset L^1(\lambda)$ and it is clear
that $L^\infty(\lambda,X) \subset E(X) \subset L^1(\lambda,X)$.
Since $X$ is weakly sequentially complete, the space $E(X)$
is weakly sequentially complete (see \cite{T2}).

Let $(f_n)_n$ be a bounded sequence in $E(X)$ that is
 equivalent to the unit vector basis of $\ell^1$.
 We will show that $(f_n)_n$ is not a (V$^*$)-set.
  If $(f_n)_n$ is not uniformly integrable then $(f_n)_n$ cannot be  
a
  (V$^*$)-set (see Proposition 3.1 of \cite{B1}) so we will assume  
that
$(f_n)_n$
  is uniformly integrable.

  By Talagrand's theorem, there exists a sequence
  $g_n \in \text{conv}(f_n, f_{n+1},\dots)$ and a measurable subset
   $\Omega^\prime$ of $\Omega$, with $\lambda(\Omega^\prime)>0$ and
   such that for each $\omega \in \Omega^\prime$, there exists
   $k \in\N$ so that $(g_n(\omega)_{n \geq k}$ is equivalent to the
   unit vector basis of $\ell^1$ in $X$. Define
   $$\varphi_n = g_n \chi_{\Omega^\prime}, \ \ \ n \in \N;$$
   Applying Theorem 1. to the sequence $(\varphi_n)_n$, there exist
   $C$ and $L$ measurable subsets
   of $\Omega$ with $\lambda(C \cup L)=1$ and a sequence
   $\psi_n \in \text{conv}(\varphi_n,\varphi_{n+1},\dots)$ so that
   \begin{itemize}
   \item[(1)] for $\omega \in C$, \ $(\psi_n(\omega) \otimes e_n)_n$  
is
 weakly Cauchy in $X \widehat{\otimes}_\pi c_0$;
\item[(2)] for $\omega \in L$, there exists $k \in \N$ so that
$(\psi_n(\omega) \otimes e_n)_{n \geq k}$ is equivalent to the unit
vector basis of $\ell^1$ in $X \widehat{\otimes}_\pi c_0$.
\end{itemize}

\noindent Case 1: Assume that $\lambda(L)>0$:

\noindent It is clear (see for instance \cite{T2}) that the sequence
$(\psi_n \otimes e_n)_n$ is equivalent to the $\ell^1$ basis in
$L^1(\lambda,X \widehat{\otimes}_\pi c_0)$ and by identification,  
the
sequence $(\psi_n \otimes e_n)_n$ is equivalent to the $\ell^1$ basis  
in
$L^1(\lambda,X) \widehat{\otimes}_\pi c_0$ so it cannot be a weakly  
null
sequence. Therefore the sequence $(\psi_n)_n$ contains a subsequence
that is equivalent to a complemented copy of $\ell^1$ in
$L^1(\lambda,X)$. Now since the inclusion map from $E(X)$ into
$L^1(\lambda,X)$ is continious, the sequence $(\psi_n)_n$ contains a
subsequence that is equivalent to a complemented copy of $\ell^1$ in
$E(X)$. As a consequence, the set
$\{\psi_n;\ n\leq 1\}$ ( and hence $\{\varphi_n;\ n \geq 1 \}$) is  
not a
(V$^*$)-set  which implies of course that the set $\{f_n, \ n \geq 1  
\}$
is not a (V$^*$)-set.

\noindent Case 2: Assume that $\lambda(L)=0$.

Since $\lambda(C \cup L)=1$ we have $\lambda(C)=1$. Note that for  
each
$\omega \in \Omega^\prime$, the sequence $(g_n(\omega))_{n \geq k}$  
is
equivalent to the unit vector basis of the $\ell^1$ for some
$k \in \N$. Now
since $\psi_n(\omega) \in
\text{conv}(g_n(\omega),g_{n+1}(\omega),\dots)$,
$(\psi_n(\omega))_{n\geq k}$ is equivalent to the unit vector basis  
of
$\ell^1$ and since $X$ has property (V$^*$), the sequence
$(\psi(\omega))_{n \geq k}$ contains a subsequence equivalent to a
complemented copy of $\ell^1$ and therefore the sequence  
$(\psi_n(\omega)
\otimes e_n)_n$ cannot be a weakly null sequence in
$X \widehat{\otimes}_\pi c_0$.
In the other hand $(\psi_n(\omega) \otimes e_n)_n$ is weakly Cauchy  
(by
the definition of $C$) so for each $\omega \in \Omega^\prime$ fixed,
there exists an operator $T \in \cal L(X,\ell^1)_1$ so that
$\lim\limits_{n \to \infty} \langle T(\psi_n(\omega)),e_n \rangle  
>0$.
Now we shall choose the operator above measurably using the  
following
proposition.
\begin{prop}
There exists a map $T: \Omega \to \cal L(X,\ell^1)_1$ with the  
following
properties:
\begin{itemize}
\item[(a)] $T(\omega)=0$\quad  $\omega \in \Omega \setminus
\Omega^\prime$;
\item[(b)] $\lim\limits_{n \to \infty} \langle
T(\omega)\left(\psi_n(\omega)\right), e_n \rangle >0 \quad \omega
\in \Omega^\prime$;
\item[(c)] The map $\omega \to T(\omega)x$ is norm-measurable for  
each $x
\in X$.
\end{itemize}
\end{prop}

We need few steps to prove the proposition.

Notice first that since $X$ is separable so is the space $X
\widehat{\otimes}_\pi c_0$
and therefore the unit ball of its dual $\cal L(X,\ell^1)_1$ is  
compact
metrizable for the weak$^*$-topology (in particular it is a Polish
space). The space $\cal L(X,\ell^1)_1 \times (X
\widehat{\otimes}_\pi c_0)^\N$ with the product topology is a Polish
space and let
$\cal A$ be a subset of
$\cal L(X,\ell^1)_1 \times (X \widehat{\otimes}_\pi
c_0)^\N$ defined as follows:
$$\{T,(\xi_n)_n \} \in \cal A
\Leftrightarrow \lim_{n \to \infty} \langle T, \xi_n \rangle >0.$$
The set $\cal A$ is trivially an Borel subset of
$\cal L(X,\ell^1)_1 \times (X \widehat{\otimes}_\pi
c_0)^\N$.
Let $\Pi: \cal L(X,\ell^1)_1 \times (X\widehat{\otimes}_\pi c_0)^\N
\to
(X\widehat{\otimes}_\pi c_0)^\N$ be  the 2nd projection; the  
operator
$\Pi$ is
of course continuous and therefore $\Pi(\cal A)$ is analytic. By  
Theorem
8.5.3 of \cite{CO}, there is a universally measurable map $\Theta:
\Pi(\cal A) \to \cal L(X,\ell^1)_1$ such that the graph of $\Theta$  
is a
subset of $\cal A$. Notice also that for $\omega \in \Omega^\prime$,  
we
have by the above argument that  the sequence
$(\psi_n(\omega) \otimes e_n)_n$ belongs to $\Pi(\cal A)$.
Now we define $T$ as follows:
$$
 T(\omega)=\begin{cases} \Theta\left( (\psi_n(\omega) \otimes
e_n)_n \right)   &\text{ for}\  \omega \in \Omega^\prime \\
                   0   &\text{ otherwise}.
\end{cases}
$$
The map $T$ is the composition of a universally measurable map
$\Theta$ and the $\lambda$- measurable map $\omega \to  
(\psi_n(\omega)
\otimes e_n)_n $ so it is measurable for the
weak$^*$-topology. Now for any $x \in X$, the map $\omega \to
T(\omega)x$ is a $\ell^1$-valued map and is weak$^*$-scalarly
measurable and since $\ell^1$ is a separable dual, it is  
norm-measurable.
Now for $\omega \in \Omega^\prime$, we get that
 $$\langle T(\omega), (\psi_n(\omega) \otimes e_n)_n
 \rangle \in \cal A $$
which by the definition of $\cal A$ is equivalent to:
$\lim_{n \to \infty}\langle T(\omega)(\psi_n(\omega)), e_n \rangle >0  
$
 and the proposition is proved. //

To finish the proof of the theorem, let  
$\gamma(\omega)=\lim\limits_{n
\to \infty} \langle T(\omega)\left(\psi_n(\omega)\right), e_n  
\rangle$.
The map $\omega \to \gamma(\omega)$ is measurable and for each  
$\omega
\in \Omega^\prime$, $\gamma(\omega)>0$. Now define $S: E(X)\to
\ell^1$ as follows:
$$S(f)= \text{Bochner}-\int_{\Omega^\prime} T(\omega)(f(\omega))\
d\lambda(\omega)$$
for each $f \in L^1(\lambda,X)$. The operator $S$ is linear and
$||S||\leq 1$ and it is easy to verify that
$$\lim_{n \to \infty} \langle S(\psi_n), e_n \rangle
=\int_{\Omega^\prime} \gamma(\omega)\ d\lambda(\omega)= \gamma >0.$$
So there exists $N \in \N$ so that for $n \geq N$, $\langle  
S(\psi_n),
e_n \rangle > \gamma/2 $ and by Proposition 1., $(\psi_n)_{n \geq N}$  
is
equivalent to a complemented copy of $\ell^1$ and therefore the set
$\{f_n;\ n\geq 1 \}$ is not a (V$^*$)-set. The proof is complete.
\end{pf}

Let us finish by asking the following question:
 Let $(\Omega,\Sigma)$ be a measure space and
$Y$ be a Banach space. We denote by $M(\Omega,Y)$ the space of
$Y$-valued countably additive measures whith bounded variation
endowed with the variation norm.

\noindent{\bf Question:} Assume that $Y=X^*$ is a dual space. Does
property (V$^*$) pass from $Y$ to $M(\Omega,Y)$?

Note that for a non dual space, the answer is negative: Talagrand
constructed in \cite{T5} a Banach lattice $E$ that does not contain
$c_0$ (so has property (V$^*$) by \cite{SS4}) but $M(\Omega,E)$
contains $c_0$ (hence failing property (V$^*$)).

\noindent {\bf Acknowledgement}: This work was completed while the
author was
visiting the Bowling Green State University, OH. The author wishes
 to thank Professor N. Carothers for his hospitality.

\bibliographystyle{plain}
\bibliography{narciref}

\end{document}